\documentclass[12pt]{article}
\usepackage{graphicx}
\usepackage{amsmath,amsthm,amssymb,enumerate}
\usepackage{euscript,mathrsfs}
\usepackage{color}
\usepackage{dsfont}
\usepackage{url}
\usepackage[left=2cm,right=2cm,top=3.5cm,bottom=3.5cm]{geometry}
\usepackage{color}
\usepackage[framemethod=tikz]{mdframed}
\allowdisplaybreaks
\usepackage{esint}

\usepackage{soul}

\catcode`\@=11 \@addtoreset{equation}{section}

\catcode`\@=12

\newtheorem{Theorem}{Theorem}[section]
\newtheorem{Proposition}[Theorem]{Proposition}
\newtheorem{Lemma}[Theorem]{Lemma}
\newtheorem{Corollary}[Theorem]{Corollary}

\theoremstyle{definition}
\newtheorem{Definition}[Theorem]{Definition}

\newtheorem{Remark}[Theorem]{Remark}

\newcommand{\bTheorem}[1]{
	\begin{Theorem} \label{T#1} }
	\newcommand{\eT}{\end{Theorem}}

\newcommand{\bProposition}[1]{
	\begin{Proposition} \label{P#1}}
	\newcommand{\eP}{\end{Proposition}}

\newcommand{\bLemma}[1]{
	\begin{Lemma} \label{L#1} }
	\newcommand{\eL}{\end{Lemma}}

\newcommand{\bCorollary}[1]{
	\begin{Corollary} \label{C#1} }
	\newcommand{\eC}{\end{Corollary}}

\newcommand{\bRemark}[1]{
	\begin{Remark} \label{R#1} }
	\newcommand{\eR}{\end{Remark}}

\newcommand{\bDefinition}[1]{
	\begin{Definition} \label{D#1} }
	\newcommand{\eD}{\end{Definition}}

\newcommand{\Ds}{\mathbb{D}_x}

\newcommand{\avintO}[1]{\fint_{\Omega} #1 \dx}

\newcommand{\bfphi}{\boldsymbol{\varphi}}

\newcommand{\bFormula}[1]{
	\begin{equation} \label{#1}}
	\newcommand{\eF}{\end{equation}}

\newcommand{\Ov}[1]{\overline{#1}}

\newcommand{\vr}{\varrho}

\newcommand{\tvr}{\wtilde \vr}
\newcommand{\tvu}{{\wtilde \vu}}

\newcommand{\vu}{\vc{u}}
\newcommand{\vm}{\vc{m}}

\newcommand{\vc}[1]{{\bf #1}}

\newcommand{\Div}{{\rm div}_x}
\newcommand{\Grad}{\nabla_x}

\newcommand{\dx}{\,{\rm d} {x}}

\newcommand{\dt}{\,{\rm d} t }

\newcommand{\vU}{\vc{U}}

\newcommand{\intO}[1]{\int_{\Omega} #1 \ \dx}

\newcommand{\intTd}[1]{\int_{\mathbb{T}^2} #1 \ \dx}

\newcommand{\D}{{\rm d}}

\newcommand{\ep}{\varepsilon}

\newcommand{\I}{\mathbb{I}}

\newcommand{\br}{ \nonumber \\ }

\def\softd{{\leavevmode\setbox1=\hbox{d}%
		\hbox to 1.05\wd1{d\kern-0.4ex{\char039}\hss}}}
\definecolor{Cgrey}{rgb}{0.85,0.85,0.85}
\definecolor{Cblue}{rgb}{0.50,0.85,0.85}
\definecolor{Cred}{rgb}{1,0,0}
\definecolor{fancy}{rgb}{0.10,0.85,0.10}
\definecolor{amaranth}{rgb}{0.9, 0.17, 0.31}

\newcommand\Cbox[2]{%
	\newbox\contentbox%
	\newbox\bkgdbox%
	\setbox\contentbox\hbox to \hsize{%
		\vtop{
			\kern\columnsep
			\hbox to \hsize{%
				\kern\columnsep%
				\advance\hsize by -2\columnsep%
				\setlength{\textwidth}{\hsize}%
				\vbox{
					\parskip=\baselineskip
					\parindent=0bp
					#2
				}%
				\kern\columnsep%
			}%
			\kern\columnsep%
		}%
	}%
	\setbox\bkgdbox\vbox{
		\color{#1}
		\hrule width  \wd\contentbox %
		height \ht\contentbox %
		depth  \dp\contentbox
		\color{black}
	}%
	\wd\bkgdbox=0bp%
	\vbox{\hbox to \hsize{\box\bkgdbox\box\contentbox}}%
	\vskip\baselineskip%
}

\mdfdefinestyle{MyFrame}{%
	linecolor=black,
	outerlinewidth=1pt,
	roundcorner=5pt,
	innertopmargin=\baselineskip,
	innerbottommargin=\baselineskip,
	innerrightmargin=10pt,
	innerleftmargin=10pt,
	backgroundcolor=white!20!white}

\newcommand{\wtilde}{\widetilde}

\allowdisplaybreaks

\begin{document}


\title{\bf Data assimilation for the barotropic Navier--Stokes system}

\author{Eduard Feireisl 	\thanks{
		The work of E.F.\ was partially supported by the
		Czech Sciences Foundation (GA\v CR), Grant Agreement
		24--11034S. The Institute of Mathematics of the Academy of Sciences of
		the Czech Republic is supported by RVO:67985840.
		E.F.\ is a member of the Ne\v cas Center for Mathematical Modelling.} }

\date{}

\maketitle

\centerline{Institute of Mathematics of the Academy of Sciences of the Czech Republic}
\centerline{\v Zitn\' a 25, CZ-115 67 Praha 1, Czech Republic}
\centerline{feireisl@math.cas.cz}

\begin{abstract}
	
We consider a continuous data assimilation method for the barotropic Navier--Stokes system. The observed solution is supposed to be bounded 
on the whole time period of observation, while the synchronized solution, usually provided by a numerical method, belongs to the class of dissipative solutions 
that is considerably larger than the class of conventional weak solutions. A complete synchronization is shown on any compact prediction interval provided 
the nudging parameters are chosen appropriately.	
	
\end{abstract}	

\section{Introduction}
\label{i}

Recently, there has been a considerable effort to establish 
a mathematical foundation of the continuous/discrete data assimilation 
techniques widely used in the real world applications of fluid mechanics, in particular in meteorology. Initiated by the seminal work of Azouani, Olson, and Titi \cite{AzOlTi}, a substantial number of results has appeared concerning mostly the 2D incompressible Navier--Stokes system, see Biswas et al \cite{BiFoMoTit}, Foias 
Mondaini and Titi, Foias, and Mondaini \cite{FoMoTi}, 
or the related Oberbeck--Boussinesq approximation examined by 
Farhat, Jolly, and Titi \cite{FaJoTi}. Applying similar techniques to the physically relevant 3D geometry is hampered by the well known difficulties concerning the well posedness of the Navier--Stokes system. 
Balakrishna and Biswas \cite{BalBis2}, \cite{BalBis} proposed a new regularity criterion for the 3D incompressible Navier--Stokes system based on the knowledge of the interpolants of the 
observed solution that can remedy the estimates in the data assimilation period, however, this piece of information is still missing in the subsequent prediction (forecast) period. 

The aim of the data assimilation process is to recover the \emph{observed solution} of a physically grounded mathematical model through the knowledge of the observation data collected in the \emph{assimilation period}. Specifically, the observed solution exists in the \emph{observation period} - the time interval $(T^-, T^+)$, and the observable data are available in the data assimilation (pre-forecast) period $[T_0, T]$, $T^- < T_0 < T < T^+$. The  goal is to construct a \emph{synchronized solution} defined in the time interval $[T_0, T^+)$ that provides the desired approximation of the observed solution in the \emph{prediction} (forecast) period $(T, T^+)$. To simplify notation, we set the beginning 
of the data assimilation period $T_0 = 0$.

The known results are based on the existence of a finite number of determining modes available for \emph{dissipative} systems, see Bercovici et al. \cite{BeCoFoMa} , Constantin et al. \cite{CFT}, \cite{CoFoMaTe}. The synchronized problem accommodates the data through 
nudging terms that force the observed and synchronized solution to approach each other exponentially fast in the data assimilation period, where it is customary to set $T = \infty$. In practice, however, 
the data assimilation (pre-forecast) period $(0,T)$ is quite short, in particular much shorter than the prediction (forecast) period 
$(T, T^+)$, see Lean et al. \cite{LBHBM}. Consequently, to ensure reliable forecast initial data, the proximity of the observed and synchronized solution must be reached at a finite time $T > 0$ by careful tuning of the nudging parameters. 

Our objective is to show convergence (synchronization) of the data assimilation method for the compressible Navier--Stokes system. Although thermal effects are completely ignored in the present study, we believe that the same approach can be applied to the full (temperature dependent) Navier--Stokes--Fourier system physically relevant in meteorological models. Our results are conditioned by two working hypotheses imposed on the observed solution: 
\begin{itemize}
\item The observed solution is  bounded in the whole observation period $(T^-,T^+)$. 
\item The observed solution is regular (smooth) at the initial time $T^-$.
\end{itemize}	
At least from the practical point of view, the above hypotheses seem quite realistic. 

Our strategy is to use first the conditional regularity results  available for the compressible Navier--Stokes system to conclude that the above hypotheses impose bounds on the derivatives of the observed solution, in particular, it is a strong solution in the whole time interval $(T^-,T^+)$. The second step is to use the relative energy inequality to estimate the distance between the observed and synchronized solution in the interval $(0, T^+)$. 

In contrast with the observed solution, the synchronized solution,
obtained in practice by means of a numerical simulation, will be understood in a general sense that goes beyond the standard framework of weak solutions. We consider
the so-called dissipative solution of the Navier--Stokes system augmented by 
nudging terms in the data assimilation period $(0,T)$. 
This type of solution, introduced in \cite{AbbFeiNov} or \cite[Chapter 5]{FeLMMiSh}, is a natural limit of consistent numerical approximations.  
The theory of dissipative solutions was developed in \cite[Chapter 5]{FeLMMiSh} with the principal aim to extend the celebrated \emph{Lax Equivalence Theorem} \cite{LaxR}
\[
\mbox{stability} + \mbox{consistency}
\ \Leftrightarrow \ \mbox{convergence}
\]
to a class 
of nonlinear systems arising in fluid dynamics.

The desired proximity of the two solution in the prediction period $(T, T^+)$ is achieved by careful tuning the nudging terms. We consider quite realistic data approximation based on collecting the exact values of the observed solution at a specific points of the space-time, in particular the network of the observation points may be time--dependent. 

The paper is organized as follows. In Section \ref{M}, we formulate the problem and state our main results. In Section \ref{P}, 
we recall the basic properties of strong and weak solutions to the barotropic Navier--Stokes system, in particular the conditional regularity estimates and the relative energy inequality. Section \ref{b}
is the heart of the paper, here we establish the estimates on the distance between the observed and synchronized solution in the 
data assimilation period. Finally, in Section \ref{C}, we show the observed and synchronized solution remain close in the prediction period $(T, T^+)$. 

\section{Problem formulation and the main result}
\label{M}

The barotropic \emph{Navier-Stokes system} governs the time evolution of the mass density $r = r(t,x)$ and the velocity $\vU = \vU(t,x)$ a compressible viscous fluid: 
\begin{align} 
	\partial_t r + \Div (r \vU) &= 0, \br
	\partial_t ( r \vU) + \Div (r \vU \otimes \vU) +  \Grad p(r) &= \Div \mathbb{S}(\Ds \vU) + r \vc{g}, \br
	\mathbb{S} &= \mu \left( \Grad \vU + \Grad^t \vU - \frac{2}{3} \Div \vU \mathbb{I} \right) + \lambda \Div \vU \mathbb{I} ,\ \mu > 0,\ \mu + \lambda \geq 0, 
	\label{O3} 	
\end{align}	
confined to a bounded spatial domain $\Omega \subset R^3$. For the sake of simplicity, we consider the no-slip boundary conditions for the velocity 
\begin{equation} \label{O4}
	\vU|_{\partial \Omega} = 0, 	
\end{equation}	
noting that our approach can be adapted to the more realistic inflow--outflow boundary conditions in a straightforward manner, see 
the monograph \cite{FeiNovOpen}. The motion is activated by an unspecified driving force $\vc{g} = \vc{g}(t,x)$ representing the effect of the ``outer world''. 

It is known that the associated initial-boundary value problem admits a global in time \emph{weak} solution as long as the pressure $p(r) \approx r^\gamma$, with $\gamma > \frac{3}{2}$, 
see \cite{EF70}, Lions \cite{LI4}. We consider the observed solution $(r,\vU)$ solving the problem \eqref{O3}, \eqref{O4} in the observation period $(T^-, T^+)$.
\begin{equation} \label{O1}
	- \infty < T^- < 0 < T^+ < \infty.
\end{equation}

\subsection{Data sampling -- interpolation operators}  
\label{d}

Fixing the data assimilation period $[0,T]$, $0 < T < T^+$, we introduce the interpolation operators $I_\delta$. For each $\delta > 0$, we consider a decomposition 
of the space--time cylinder:
\begin{align}
	[0,T] \times \Ov{\Omega} = \cup_{i = 1}^{N(\delta)} \Ov{Q}_{i, \delta}, \ Q_{i,\delta} \cap Q_{j, \delta} = \emptyset \ \mbox{for}\ i \ne j, 
	\ {\rm diam}[ Q_{i, \delta} ] \leq \delta,
	\label{d1}
\end{align}
where $Q_{i, \delta}$ are domains, 
together with the control points 
\begin{equation} \label{d2}
	y_{i,\delta} = (t_{i, \delta}, x_{i, \delta}) \in Q_{i, \delta},\ 
	i = 1, \dots, N(\delta).
\end{equation}	

Given a continuous function $\theta \in C([0,T] \times \Ov{\Omega})$, we set 
\begin{equation} \label{d3}
	I_\delta [\theta] = \sum_{i = 1}^{N(\delta)} \theta (y_{i, \delta}) 
	\mathds{1}_{Q_{i, \delta}} \in L^\infty((0,T) \times \Omega).	
\end{equation}

\subsection{Synchronized system - data assimilation}
\label{S}

The synchronized system is defined in the time interval $(0, T^+)$, where the time interval $(T, T^+)$ is called prediction period. 
It is the barotropic Navier--Stokes system \eqref{O3} perturbed by nudging terms active in the data assimilation period $(0,T)$: 
\begin{align} 
	\partial_t \vr + \Div (\vr \vu) &= - \Lambda_\vr (\vr - I_\delta[r]) \mathds{1}_{(0,T)}, \br 
	\partial_t (\vr \vu) + \Div (\vr \vu \otimes \vu) + \Grad p(\vr) &= 
	\Div \mathbb{S}(\Ds \vu) + \vr \vc{g} - \Lambda_\vu (1 + \vr) (\vu - I_\delta[\vU] )  \mathds{1}_{(0,T)} 
	\label{S1}
\end{align}	
in $(0, T^+) \times \Omega$, with the no--slip boundary conditions 
\begin{equation} \label{S2}
	\vu|_{\partial \Omega} = 0. 	 	 
\end{equation}

The initial conditions for the synchronized system can be taken arbitrarily. For definiteness and also simplicity, we consider 
\begin{equation} \label{S3}
	\vr(0, \cdot) = \vr_0 = \avintO{ r(T^-, \cdot) }, \ \vm_0 = \vr \vu(0, \cdot) = 0
\end{equation}	
bearing in mind that the total mass of the observed solution, 
\[
M_0 = \intO{ r(t, \cdot) },
\]
is a constant of motion. 
Here $\Lambda_\vr, \Lambda_\vu > 0$ are a positive relaxation parameters to be fixed below.

Unlike the ``standard'' nudging $I_\delta [\vr] - I_\delta [r]$, 
$I_\delta[\vu] - \I_\delta[\vU]$, the interpolation operators in \eqref{S1} act only on the observed solution. There are essentially two reasons for this strategy. The first one, rather obvious, is the low regularity of the synchronized solution for the interpolation $I_\delta$ to be well defined. The second, more subtle, is to compensate the lack of dissipation provided by the physical model in particular at the level of the mass density.

We consider the so-called dissipative solutions of the synchronized problem 
proposed in \cite{AbbFeiNov}. First, let us introduce the total energy associated to the compressible Navier--Stokes system,
\begin{equation} \label{S4}
	E(\vr, \vm) = \left\{ \begin{array}{l} \frac{1}{2} \frac{|\vm|^2}{\vr} + P(\vr) \ \mbox{if}\ \vr > 0, \\ \\
		0 \ \mbox{if}\ \vr = 0, \vm = 0,\\ \\ \infty \ \mbox{otherwise} \end{array} \right.,\ \mbox{where}\ 
	P'(\vr) \vr - P(\vr) = p(\vr),\ \vm = \vr \vu.
\end{equation} 
It is easy to check that $E: R^4 \to [0, \infty)$ is a convex l.s.c function as soon as 
\begin{equation} \label{S4a}
	p'(\vr) > 0 \ \mbox{for all}\ \vr > 0.
\end{equation}	
In addition, we require 
\begin{equation} \label{S4A} 
p(0) = P(0) = 0,\ \vr \mapsto p(\vr)\ \mbox{convex},\ 
\vr \mapsto P(\vr) - a p(\vr) \ \mbox{convex for a positive constant}\ a > 0.
\end{equation}
It is easy to check that hypothesis \eqref{S4A} is satisfied if $p \in C^2(0, \infty)$, and 
\[
p'(\vr) \geq a p''(\vr) \vr > 0 \ \mbox{for all}\ \vr > 0.
\]

\begin{Definition}[\bf Dissipative solution] \label{SD1}

We say the $(\vr, \vu)$ is \emph{dissipative solution} to the synchronized problem \eqref{S1}--\eqref{S3} in 
$(0,T^+) \times \Omega$ if the following holds true:	
\begin{itemize}
\item {\bf Regularity.} The density $\vr \in C_{\rm weak}([0, T^+]; L^\gamma (\Omega))$ for a certain $\gamma > 1$, the velocity $\vu \in L^2(0, T^+; W^{1,2}_0(\Omega; R^3))$, the 
momentum $\vm = \vr \vu \in C_{\rm weak}([0, T^+]; L^{\frac{2 \gamma}{\gamma + 1}}(\Omega; R^3))$, the total energy 
$E(\vr, \vm)$ belongs to $L^\infty(0, T^+; L^1(\Omega))$. 

\item {\bf Synchronized equation of continuity.} The integral identity 
\begin{equation} \label{SS1}
\int_0^{T^+} \intO{ \Big[ \vr \partial_t \varphi + \vr \vu \cdot \Grad \varphi \Big] } \dt = 
\Lambda_\vr \int_0^T \intO{ (\vr - I_\delta[r]) \varphi } \dt - \intO{ \vr_0 \varphi (0, \cdot)} 
\end{equation}	
holds for any $\varphi \in C^1_c([0, T^+) \times \Ov{\Omega})$. 

\item {\bf Synchronized momentum equation.} The integral identity 	
\begin{align} 
\int_0^{T^+} &\intO{ \Big[ \vr \vu \cdot \partial_t \bfphi + \vr (\vu \otimes \vu): \Grad \bfphi + 
p(\vr) \Div \bfphi \Big] } \dt \br &= \int_0^{T^+} \intO{ \Big[  \mathbb{S}(\Ds \vu) : \Grad \bfphi - \vr \vc{g} \cdot \bfphi \Big] } \dt 	
- \int_0^{T^+} \left( \int_{\Omega} \Grad \bfphi : \D \mathcal{R} \right) \dt \br
&+ \Lambda_\vu \int_0^T \intO{ (1 + \vr) (\vu - I_\delta [\vU] ) \cdot \bfphi } \dt - \intO{ \vm_0 \cdot \bfphi (0, \cdot) } 	
\label{SS2}
\end{align}	
holds for any $\bfphi \in C^1_c([0, T^+) \times \Omega)$, where $\mathcal{R} \in L^\infty(0,T; \mathcal{M}^+ (\Ov{\Omega}; R^3))$.

\item {\bf Energy inequality.}
The integral inequality 
\begin{align} 
- &\int_0^{T^+} \partial_t \psi  \left[ \intO{ E(\vr, \vr \vu) } \dt + a \int_{\Ov{\Omega}} \D\ {\rm trace}[\mathcal{R}] \right] + \int_0^{T^+} \psi \intO{ \mathbb{S}(\Ds \vu) : 
		\Ds \vu } \dt\br  &+ \Lambda_\vu \int_0^T \psi  \intO{|\vu|^2 } \dt + (\Lambda_\vu - \Lambda_\vr) 
	\int_0^T \psi \intO{ \vr |\vu|^2 } \dt + \frac{1}{2} \Lambda_\vr \int_0^T \psi \intO{ I_\delta[r] |\vu|^2 } \dt \br &+ \frac{1}{2} \Lambda_\vr \int_0^T \psi \intO{ \vr |\vu|^2} \dt
	  + \Lambda_\vr \int_0^T \psi \intO{ \Big( P(\vr) - P( I_\delta[r] ) \Big) } \dt \br 
&+ a \Lambda_r \int_0^T \int_{\Ov{\Omega}} \D \ {\rm trace}[\mathcal{R}]  \dt \br&\leq 
	\int_0^{T^+} \psi \intO{ \vr \vc{g} \cdot \vu } \dt  +  \Lambda_\vu \int_0^T \psi \intO{ (1 + \vr) I_\delta [\vU] \cdot \vu }\dt + \psi(0) \intO{ E(\vr_0, \vm_0) }.
	\label{SS3}
\end{align}	
holds for any $\psi \in C^1_c[0, T^+)$, $\psi \geq 0$.
\end{itemize}
	
\end{Definition}

The specific form of the energy inequality \eqref{SS3} deserves some comments. Formally, the 
energy balance is derived by multiplying the momentum equation in \eqref{S1} on $\vu$ and using the equation of continuity: 
\begin{align}
 \frac{\D }{\dt } &\intO{ E(\vr, \vr \vu) } + \intO{ \mathbb{S}(\Ds \vu) : \Ds \vu }\br   &+ \mathds{1}_{(0,T)} \Lambda_\vr \intO{ P'(\vr) (\vr - I_\delta[r]) }  
 + \mathds{1}_{(0,T)}\frac{\Lambda_\vr}{2} \intO{ I_\delta[r] |\vu|^2 } 
  \br &+ \mathds{1}_{(0,T)} \Lambda_\vu \intO{  |\vu|^2 } + 
\mathds{1}_{(0,T)}  (\Lambda_\vu - \Lambda_\vr)\intO{ \vr  |\vu|^2 }  + \frac{\Lambda_\vr}{2} \mathds{1}_{(0,T)} \intO{ \vr  |\vu|^2 } \br &=   \intO{ \vr \vc{g} \cdot \vu } 
+ \mathds{1}_{(0,T)} \Lambda_\vu \intO{ (1 + \vr) I_\delta[\vU] \cdot \vu }. 
\nonumber	
\end{align}	
Moreover, by virtue of Fenchel--Young inequality, 
\[
P'(\vr) \vr - P'(\vr) I_\delta [r] = P(\vr) + P^* (P'(\vr)) - P'(\vr) I_\delta [r] \geq P(\vr) - P(I_\delta[r]).
\]	
where $P^*$ denotes the conjugate of $P$. Consequently, we get 
\begin{align}
	\frac{\D }{\dt } &\intO{ E(\vr, \vr \vu) } + \intO{ \mathbb{S}(\Ds \vu) : \Ds \vu }  
	+ \mathds{1}_{(0,T)}\frac{\Lambda_\vr}{2} \intO{ I_\delta[r] |\vu|^2 } 
	\br &+ \mathds{1}_{(0,T)} \Lambda_\vu \intO{  |\vu|^2 } + 
	\mathds{1}_{(0,T)}  (\Lambda_\vu - \Lambda_\vr)\intO{ \vr  |\vu|^2 } \br &+ \mathds{1}_{(0,T)} \frac{\Lambda_\vr}{2} \intO{ \vr  |\vu|^2 } 
	+ \mathds{1}_{(0,T)} \Lambda_\vr \intO{ \Big( P(\vr) - P(I_\delta[r]) \Big) }
	\br &\leq   \intO{ \vr \vc{g} \cdot \vu } 
	+ \mathds{1}_{(0,T)} \Lambda_\vu \intO{ (1 + \vr) I_\delta[\vU] \cdot \vu },
	\label{pom}
\end{align}	
which is identical with \eqref{SS3} except the integrals containing ${\rm trace}[\mathcal{R}]$. 

The concept of dissipative solution anticipates possible oscillations/concentrations that may arise in a sequence of numerical solutions. 
They are encoded in the quantity $\mathcal{R}$ usually termed Reynolds stress, 
\[
\mathcal{R} = \mbox{weak - (*) limit } \Big( \frac{\vm_n \otimes \vm_n}{\vr_n} + p(\vr_n) \Big) - \Big( \frac{\vm \otimes \vm}{\vr} + p(\vr) \Big) 
\]
for any approximate sequence $\vm_n = \vr_n \vu_n$, $\vr_n$ converging weakly to $\vm$, $\vr$, respectively.
The same effect is observed for the total energy, where 
\[
\mathcal{C} = \mbox{weak - (*) limit } \Big( \frac{1}{2} \frac{|\vm_n|^2}{\vr_n} + P(\vr_n) \Big) - \Big( \frac{1}{2} \frac{|\vm|^2}{\vr} + P(\vr) \Big).
\] 
Without loss of generality, we may assume the constant $a$ in hypothesis \eqref{S4A} satisfies $0 < a < \frac{1}{2}$. It follows 
\[
\mathcal{C} \geq a\ {\rm trace} [\mathcal{R}].
\]
Finally, we may modify the Reynolds stress by adding a spatially homogenous function $\mathcal{R}(t, \cdot) + \chi(t) \mathbb{I}$ obtaining 
\[
\int_{\Ov{\Omega}} \D \mathcal{C} = a \int_{\Ov{\Omega}} \D \ {\rm trace}[\mathcal{R}].
\]
Incorporating the oscillation defect into \eqref{pom} we get \eqref{SS3}.

Obviously, any dissipative solution is a weak solution of the synchronized system \eqref{S1}--\eqref{S3} as soon as $\mathcal{R} = 0$. The existence of 
global in time \emph{weak} solutions can be shown by the method developed in \cite{EF70}, \cite{LI4} as soon as the pressure satisfies an addition growth condition 
\begin{equation} \label{SS3a}
p(\vr) \geq \underline{p} \vr^\gamma > 0 \ \mbox{for all}\ \vr \geq 1 
\ \mbox{and}\ \gamma > \frac{3}{2}.
\end{equation}	
The dissipative solutions exist globally in time for the full range $\gamma > 1$, cf. \cite[Chapter 5]{FeLMMiSh}. Moreover, as shown in 
\cite[Chapter 7]{FeLMMiSh}, they arise as weak limits of certain consistent numerical approximations.

\subsection{Main result}
\label{m}

We consider the relative energy
\begin{equation} \label{S4b}
E \left( \vr, \vu \Big| \tvr, \tvu \right) = 
\frac{1}{2} \vr |\vu - \tvu|^2 + P(\vr) - P'(\tvr)(\vr - \tvr) - P(\tvr). 
\end{equation}	
As shown in \cite{FeJiNo}, see also a more elaborated treatment in 
\cite[Chapter 4]{FeLMMiSh}, the relative energy can be interpreted as a generalized distance - the so-called Bregmann divergence - between the quantities $(\vr, \vu)$ and $(\tvr, \tvu)$.

Before stating our main result, we introduce the ``data'' of the problem. We suppose the observed solution $(r, \vU)$ is bounded on the whole time interval $(T^-, T^+)$, specifically 
\begin{equation} \label{O2}
	0 \leq r(t,x) \leq \Ov{r},\ 
	|\vU(t,x) | \leq \Ov{U} \ \mbox{for}\ x \in \Ov{\Omega} 
	\ \mbox{uniformly for}\ T^- < t < T^+. 
\end{equation} 
In addition, we assume 
\begin{equation} \label{O6}	
\| r(T^-, \cdot) \|_{W^{1,q}(\Omega)} < \infty,\ 
r(T^-, \cdot) > 0 \ \mbox{in}\ \Ov{\Omega},\  
\| \vU(T^-, \cdot) \|_{W^{2(1 - \frac{1}{q}),q}(\Omega)} < \infty   
\end{equation}
for a certain $q > 3$ specified below, 
and set 
\[
D_0 = \inf_{t \in [T^-, 0]} \left( \| r(t, \cdot) \|_{W^{1,q}(\Omega)} + \| r^{-1}(t, \cdot) \|_{L^\infty(\Omega)} + \| \vU(t, \cdot) \|_{W^{2(1 - \frac{1}{q}),q}(\Omega)}  \right) < \infty.
\]

Finally, we denote 
\[
\Ov{g} = {\rm ess} \sup_{ (t,x) (T^-, T^+) \times \Omega} | 
\vc{g} (t,x)|, 
\]
and define 
\[
\| {\rm data} \| = D_0 + \Ov{r} + \Ov{U} + \Ov{g} + (T^+ - T^-) .
\]

We are ready to state our main result. 

\begin{Theorem}[\bf Synchronization for the compressible Navier-Stokes system] \label{mT1}
	
Let $\Omega \subset R^3$ be a bounded domain of class $C^4$. Let the pressure--density equation of state $p=p(\vr)$ satisfy 
the hypotheses \eqref{S4a}, \eqref{S4A}, together with the growth restriction \eqref{SS3a}, where $\gamma > \frac{6}{5}$.
Suppose the observed solution $(r, \vU)$ belongs to the regularity class 
\eqref{O2}, \eqref{O6}, with 
\[
q > \max \left\{ 12; \frac{6 \gamma}{5 \gamma - 6} \right\}.
\]
 Let  
$(\vr, \vu)$ be a dissipative solution of the synchronized system \eqref{S1}--\eqref{S3} in the sense of Definition \ref{SD1}.

Then there is a non--decreasing function $\Gamma: (0, \infty) \to [1, \infty)$ such that 
\begin{equation} \label{sy1}
{\rm ess} \sup_{t \in [T, T^+]} \left[ \intO{ E \left( \vr, \vu \Big| r, \vU \right) (t, \cdot) } + a \int_{\Ov{\Omega}} \D {\rm trace}[\mathcal{R}](t) \right] < \ep 	
\end{equation}	
whenever  
\begin{align} 
\left[ \frac{1}{\Lambda_\vr} + \exp \left( - {\Lambda_\vr} T \right) \right] \Gamma (\| {\rm data} \|) &< \ep, \br
\Lambda_\vu &\geq \Gamma (\| {\rm data} \|) \Lambda_\vr,  \br
\delta \Lambda_\vu \Gamma (\| {\rm data} \|) &\leq 1.
\label{sy2}	
	\end{align}	 
	
\end{Theorem}

\begin{Remark} \label{mR2}
	
The restriction $\gamma > \frac{6}{5}$ can be replaced by $\gamma > 1$ if the values of the observed solution $r(T^-, \cdot)$, $\vU(T^-, \cdot)$ are more regular.
This point is clearly indicated in the proof below, cf. Remarks \ref{IR1}, \ref{IR2}.

\end{Remark}

The rest of the paper is devoted to the proof of Theorem \ref{mT1}. We use the symbol $\Gamma$ to denote a generic non--decreasing function ranging in 
$[1, \infty)$.

\section{Preliminary results concerning the observed and synchronized system}
\label{P}

We recall some known facts concerning the properties of the observed and synchronized system.

\subsection{Conditional regularity for the observed system}

In their seminal work, Sun, Wang, and Zhang \cite{SuWaZh1} showed that any weak solution of the barotropic Navier--Stokes system emanating from the 
smooth initial data and with uniformly bounded density is in fact regular. Here, we report an $L^q$ variant of their result proved in \cite[Theorem 7.1]{AbBaChFe}.

\begin{Proposition}[\bf Conditional regularity] \label{PO1}
	For any $3 < q < \infty$, the observed solution $(r, \vU)$ 
	belongs to the class 
\begin{align} 
r &\in C([T^-,T^+]; W^{1,q}(\Omega)),\ 
\partial_t \vr \in C([T^-, T^+]; L^q(\Omega)),\ \inf_{[T^-, T^+] \times \Ov{\Omega}} r > 0, \br
\vU &\in L^q(T^-,T^+; W^{2,q} \cap W^{1,q}_0 (\Omega; R^3)), \ \partial_t \vU \in L^q(T^-; T^+; L^q(\Omega; R^3)), \br	
\label{cla1}	
\end{align}		
and the following estimate holds true:
	\begin{align} 
	\sup_{t \in [T^-,T^+]} & \| r(t, \cdot) \|_{W^{1,q}(\Omega)} 
	+ \sup_{t \in [T^-, T^+]} \| \partial_t r(t, \cdot) \|_{ L^q(\Omega)} + 
	\| r^{-1} \|_{L^\infty((T^-, T^+) \times \Omega)} \br	
	&+ \int_{T^-}^{T^+} \| \vU(t, \cdot) \|^q_{W^{2,q}(\Omega; R^3)} \dt  + 
	\int_{T^-}^{T^+} \| \partial_t \vU(t,\cdot) \|^q_{L^q(\Omega; R^3)} \dt \br 
&\quad \leq \Gamma \Big( \| {\rm data} \| \Big), 
		\label{O8}
	\end{align}	
	where $\Gamma$ is a non--decreasing function. 	
	
\end{Proposition} 	

Strictly speaking, the result proved in \cite{AbBaChFe} applies to the full Navier--Stoke--Fourier system. The present version for the barotropic case is a straightforward consequence.

\subsection{Relative energy inequality for the synchronized system} 

The relative energy inequality can be deduced in the way proposed in \cite{FeJiNo}. 
The main idea is to deduce a formula governing the time evolution 
of the quantity 
\[
E \left( \vr, \vu \Big| r, \vU \right) = \frac{1}{2} 
\vr |\vu - \vU|^2 + P(\vr) - P'(r)(\vr - r) - P(r) 
\]
in the time interval $[0, T^+]$. This can be achieved by considering $\vU$ as a test function in the weak 
formulation of the momentum balance $\eqref{SS2}$, and the quantity $\frac{1}{2}|\vU|^2 - P'(r)$ 
as a test function in the weak formulation of the equation of continuity $\eqref{SS1}$, see \cite{FeJiNo} for details. 
After a tedious but straightforward manipulation using 
the fact that $(r, \vU)$ is a strong solution of the Navier--Stokes system, we obtain 
\begin{align} 
	&\frac{\D }{\dt} \left[ \intO{ E \left(\vr, \vu \ \Big|\  r, \vU \right)} + a \int_{\Ov{\Omega}} \D \ {\rm trace}[ \mathcal{R}] \right]  + 
	\intO{ \left( \mathbb{S}(\Ds \vu) - \mathbb{S} (\Ds \vU) \right): \left( \Ds \vu - \Ds \vU \right) }\br &+ \mathds{1}_{(0,T)} \Lambda_\vu \intO{ (1 + \vr) (\vu - I_\delta [\vU] ) \cdot (\vu - \vU) } \br	&+ \mathds{1}_{(0,T)} \Lambda_\vr \intO{  \Big( P(\vr) - P'(r) (\vr- I_\delta[r]) - P(I_\delta[r]) \Big) } \br 
	&+  \mathds{1}_{(0,T)} a \Lambda_\vr \int_{\Ov{\Omega}} \D \ {\rm trace}[\mathcal{R}] \br
	&\leq - \intO{ \vr (\vU - \vu) \cdot \Grad \vU \cdot (\vU - \vu) } \br
	&\quad -  \intO{ \Div \vU \Big( p(\vr) - p'(r) (\vr - r) -  p(r) \Big) } \br
	&\quad + \int_{\Ov{\Omega}} \Grad \vU : \D \mathcal{R} \br
	&\quad + \intO{ (\vr - r) \Big( \frac{1}{r} \Div \mathbb{S}(\Ds \vU) + \vc{g} \Big) \cdot (\vU - \vu) }                  \br
	& \quad + \frac{1}{2}  \mathds{1}_{(0,T)} \Lambda_\vr \intO{ (I_\delta[r] - \vr)   ( \vU + \vu ) 
		\cdot(\vU - \vu) } \ \mbox{in}\ \mathcal{D}'(0,T^+).
	\label{S6}
\end{align} 

At this stage, it is worth observing that 
\[
\intO{ \left( \mathbb{S}(\Ds \vu) - \mathbb{S} (\Ds \vU) \right): \left( \Ds \vu - \Ds \vU \right) } \geq \nu \| \vu - \vU \|^2_{W^{1,2}_0(\Omega; R^3)} 
\]
as a direct consequence of Korn--Poincar\' e inequality. Moreover, we have 
\begin{align}
 \Big( P(\vr) - P'(r) (\vr- I_\delta[r]) - P(I_\delta[r]) \Big) &= 
  \Big( P(\vr) - P'(r) (\vr- r) - P(r) \Big) \br &- \Big(  P(I_\delta[r]) - P'(r) (I_\delta[r] - r) - P(r) \Big).
	\label{S7a}
\end{align}

Consequently, inequality \eqref{S6} can be rewritten in the form 
\begin{align} 
	&\frac{\D }{\dt} \left[ \intO{ E \left(\vr, \vu \ \Big|\  r, \vU \right) } + a \int_{\Ov{\Omega}} 
	\D \ {\rm trace}[\mathcal{R}]    \right] \br &+ \mathds{1}_{(0,T)} \Lambda_\vr \left[ \intO{ E \left(\vr, \vu \ \Big|\  r, \vU \right) } + a \int_{\Ov{\Omega}} 
	\D \ {\rm trace}[\mathcal{R} ]   \right] \br 
	&  + \nu \| \vu - \vU \|^2_{W^{1,2}_0(\Omega; R^3)} + \mathds{1}_{(0,T)} \Lambda_\vu \|\vu - \vU\|^2_{L^2(\Omega; R^3)} + \mathds{1}_{(0,T)} \frac{\Lambda_\vu}{2} 
	\intO{ \vr |\vu - \vU|^2 }
\br 
	&\quad \leq - \intO{ \vr (\vU - \vu) \cdot \Grad \vU \cdot (\vU - \vu) } \br
	&\quad \quad -  \intO{ \Div \vU \Big( p(\vr) - p'(r) (\vr - r) -  p(r) \Big) } \br
	&\quad \quad + \int_{\Ov{\Omega}} \Grad \vU : \D \mathcal{R} \br
	&\quad \quad + \intO{ (\vr - r) \Big( \frac{1}{r} \Div \mathbb{S}(\Ds \vU) + \vc{g} \Big) \cdot (\vU - \vu) }                  \br
	&\quad \quad + \frac{1}{2}  \mathds{1}_{(0,T)} \Lambda_\vr \intO{ (I_\delta[r] - \vr)   ( \vU + \vu ) 
		\cdot(\vU - \vu) } \br  &\quad \quad+  \mathds{1}_{(0,T)} \Lambda_\vr \intO{ \Big(P(I_\delta[r]) -  P'(r) (I_\delta[r] - r) - P(r) \Big)  } \br  &\quad \quad+ \mathds{1}_{(0,T)} \Lambda_\vu \intO{ (1 + \vr) (I_\delta[\vU] - \vU) \cdot (\vu - \vU) } 
	\label{S6a}
\end{align} 
in $\mathcal{D}'(0,T^+)$ as soon as 
\begin{equation} \label{S8}
	\Lambda_\vu \geq \Lambda_\vr.
\end{equation}	 

\section{Bounds on the relative energy in the data assimilation period}
\label{b}

We focus on the data assimilation period $(0,T)$, where the relative energy 
inequality is supplemented by nudging terms. First observe that 
\begin{align}
\intO{ (I_\delta[r] - \vr)   ( \vU + \vu ) 
	\cdot(\vU - \vu) } &= \intO{ (\vr - I_\delta[r]) |\vu - \vU|^2 } + 2 
\intO{ (I_\delta[r] - \vr) \vU \cdot (\vU - \vu) } \br 
&\leq \intO{ \vr  |\vu - \vU|^2 } + 2 
\intO{ (I_\delta[r] - \vr) \vU \cdot (\vU - \vu) }.
\nonumber
\end{align}
Consequently, strengthening \eqref{S8} to 
\begin{equation} \label{S8a} 
\Lambda_\vu \geq 2 \Lambda_\vr,	
\end{equation}	
we deduce from \eqref{S6a}
\begin{align} 
&\frac{\D }{\dt} \left[ \intO{ E \left(\vr, \vu \ \Big|\  r, \vU \right) } + a \int_{\Ov{\Omega}} 
\D \ {\rm trace}[\mathcal{R}]    \right] \br &+ \Lambda_\vr \left[ \intO{ E \left(\vr, \vu \ \Big|\  r, \vU \right) } + a \int_{\Ov{\Omega}} 
\D \ {\rm trace}[\mathcal{R} ]   \right] \br
	&  + \nu \| \vu - \vU \|^2_{W^{1,2}_0(\Omega; R^3)} +  \Lambda_\vu \|\vu - \vU \|^2_{L^2(\Omega; R^3)} +  \frac{\Lambda_\vu}{4} 
	\intO{ \vr |\vu - \vU|^2 }
	\br 
	&\quad \leq - \intO{ \vr (\vU - \vu) \cdot \Grad \vU \cdot (\vU - \vu) } \br
	&\quad \quad -  \intO{ \Div \vU \Big( p(\vr) - p'(r) (\vr - r) -  p(r) \Big) } \br
	&\quad \quad + \int_{\Ov{\Omega}} \Grad \vU : \D \mathcal{R} \br
	&\quad \quad + \intO{ (\vr - r) \Big( \frac{1}{r} \Div \mathbb{S}(\Ds \vU) + \vc{g} \Big) \cdot (\vU - \vu) }                  \br
	&\quad \quad +  \Lambda_\vr \intO{ (I_\delta[r] - \vr)  \vU  
		\cdot(\vU - \vu) } \br  &\quad \quad+   \Lambda_\vr \intO{ \Big(P(I_\delta[r]) -  P'(r) (I_\delta[r] - r) - P(r) \Big)  } \br  &\quad \quad+  \Lambda_\vu \intO{ (1 + \vr) (I_\delta[\vU] - \vU) \cdot (\vu - \vU) } 
	\label{S9}
\end{align} 
for in $\mathcal{D}'(0,T)$.

\subsection{Uniform bounds depending on the data}

Thanks to hypothesis \eqref{S4a}, we have
\begin{equation} \label{S10}
0 \leq a \Big( p(\vr) - p'(r) (\vr - r) - p(r) \Big) \leq \Big( P(\vr) - P'(r) (\vr - r) - P(r) \Big). 
\end{equation}	
Consequently, 
\begin{align} 
&	- \intO{ \vr (\vU - \vu) \cdot \Grad \vU \cdot (\vU - \vu) }  -  \intO{ \Div \vU \Big( p(\vr) - p'(r) (\vr - r) -  p(r) \Big) } + 
\int_{\Omega} \Grad \vU : \D \mathcal{R} \br &\leq  
C \| \Grad \vU \|_{L^\infty(\Omega; R^{3 \times 3})} \left[ \intO{ E \left(\vr, \vu \ \Big|\  r, \vU \right) } + a \int_{\Ov{\Omega}} 
\D \ {\rm trace}[\mathcal{R}]    \right].  
\label{S10a}
\end{align}	
In view of the conditional regularity bounds \eqref{O8}, we have 
\[
\vU \in C([0, T^+]; W^{2 (1 - \frac{1}{q}), q}(\Omega; R^3)),\ 
\mbox{where}\  W^{2 (1 - \frac{1}{q}), q}(\Omega; R^3)) \hookrightarrow C^1(\Ov{\Omega}; R^3) \ \mbox{provided}\ q > 5.
\]
Thus we may rewrite \eqref{S9} as 
\begin{align} 
	&\frac{\D }{\dt} \left[ \intO{ E \left(\vr, \vu \ \Big|\  r, \vU \right) } + a \int_{\Ov{\Omega}} 
	\D \ {\rm trace}[\mathcal{R}]    \right] + \Lambda_\vr \left[ \intO{ E \left(\vr, \vu \ \Big|\  r, \vU \right) } + a \int_{\Ov{\Omega}} 
	\D \ {\rm trace}[\mathcal{R} ]   \right] \br 
	&  + \nu \| \vu - \vU \|^2_{W^{1,2}_0(\Omega; R^3)} +  \Lambda_\vu \|\vu - \vU \|^2_{L^2(\Omega; R^3)} +  \frac{\Lambda_\vu}{4} 
	\intO{ \vr |\vu - \vU|^2 }
	\br 
	&\quad \leq \Gamma (\| {\rm data} \|)\left[ \intO{ E \left(\vr, \vu \ \Big|\  r, \vU \right) } + a \int_{\Ov{\Omega}} 
	\D \ {\rm trace}[\mathcal{R}]    \right]  \br
	&\quad \quad + \intO{ (\vr - r) \Big( \frac{1}{r} \Div \mathbb{S}(\Ds \vU) + \vc{g} \Big) \cdot (\vU - \vu) }                  \br
	&\quad \quad +  \Lambda_\vr \intO{ (\vr - r)  \vU  
		\cdot(\vu - \vU) } \br&\quad \quad + \Lambda_\vr \intO{ (I_\delta[r] - r)  \vU  
		\cdot(\vU - \vu) } \br &\quad \quad+   \Lambda_\vr \intO{ \Big(P(I_\delta[r]) -  P'(r) (I_\delta[r] - r) - P(r) \Big)  } \br  &\quad \quad+  \Lambda_\vu \intO{ (1 + \vr) (I_\delta[\vU] - \vU) \cdot (\vu - \vU) }.
	\label{b6}
\end{align} 

Next, we decompose 
\begin{align}
&\intO{ (\vr - r) \Big( \frac{1}{r} \Div \mathbb{S}(\Ds \vU) + \vc{g} \Big) \cdot (\vU - \vu) }  = 
\int_{\vr \leq 2r} (\vr - r) \Big( \frac{1}{r} \Div \mathbb{S}(\Ds \vU) + \vc{g} \Big) \cdot (\vU - \vu) \dx\br &\quad+ 
\int_{\vr > 2r} (\vr - r) \Big( \frac{1}{r} \Div \mathbb{S}(\Ds \vU) + \vc{g} \Big) \cdot (\vU - \vu) \dx, 
\nonumber
\end{align}
where, by means of H\" older inequality, 
\begin{align}
\int_{\vr \leq 2r}&(\vr - r) \Big( \frac{1}{r} \Div \mathbb{S}(\Ds \vU) + \vc{g} \Big) \cdot (\vU - \vu) \dx \br &\leq c(\omega)
\left\| \frac{1}{r} \Div \mathbb{S}(\Ds \vU) + \vc{g} \right\|_{L^3(\Omega)}^2  \int_{\vr \leq 2 r} |\vr - r |^2 \dx + 
\omega \| \vu - \vU \|^2_{L^6(\Omega; R^d)} 
\label{b7}
\end{align}
for arbitrary $\omega > 0$. Thus adjusting $\omega = \omega(\nu)$ appropriately, we get 
\begin{align}
c(\omega)
&\left\| \frac{1}{r} \Div \mathbb{S}(\Ds \vU) + \vc{g} \right\|_{L^3(\Omega)}^2  \int_{\vr \leq 2 r} |\vr - r |^2 \dx + 
\omega \| \vu - \vU \|^2_{L^6(\Omega; R^d)} \br
&\leq c(\nu, \Ov{r}) \left\| \frac{1}{r} \Div \mathbb{S}(\Ds \vU) + \vc{g} \right\|_{L^3(\Omega)}^2 \intO{ E \left(\vr, \vu \Big| r, \vU \right) } + \frac{\nu}{2} \| \vu - \vU\|^2_{W^{1,2}_0(\Omega; R^d)}, 
\label{b7a}
\end{align}
where we have used the embedding $W^{1,2} \hookrightarrow L^6$.

As for the second integral, we fix $q = 12$, $r = \frac{12}{5}$, and use H\" older inequality obtaining 
\begin{align}
\int_{\vr > 2r} &(\vr - r) \Big( \frac{1}{r} \Div \mathbb{S}(\Ds \vU) + \vc{g} \Big) \cdot (\vU - \vu) \dx \leq 
\int_{\vr > 2r} \sqrt{\vr} \Big| \frac{1}{r} \Div \mathbb{S}(\Ds \vU) + \vc{g} \Big| \cdot \sqrt{\vr} (\vU - \vu) \dx\br \leq 
&\left( \int_{\vr > 2r} \vr^{\frac{6}{5}} \dx \right)^{\frac{5}{12}} 
\left\| \frac{1}{r} \Div \mathbb{S}(\Ds \vU) + \vc{g} \right\|_{L^q(\Omega; R^d)} \left\|  \sqrt{\vr} (\vU - \vu)  \right\|_{L^2(\Omega); R^d)} \br 
&\leq \left\| \frac{1}{r} \Div \mathbb{S}(\Ds \vU) + \vc{g} \right\|_{L^9(\Omega; R^d)}^2 \intO{ E\left( \vr, \vu \Big| r, \vU \right) } + \left(  \intO{ E\left( \vr, \vu \Big| r, \vU \right) } \right)^{\frac{5}{6}}.
\label{b8}
\end{align}	

\begin{Remark} \label{IR1}
Note this argument can be carried out for any $\gamma > 1$,                                                                                                                                                                                                                                                                                                    with $r = 2 \gamma$ and $q = \frac{\gamma - 1}{2 \gamma}$.
\end{Remark}	

Summing up \eqref{b7}--\eqref{b8}, we may rewrite \eqref{b6} in the form
\begin{align} 
	&\frac{\D }{\dt} \left[ \intO{ E \left(\vr, \vu \ \Big|\  r, \vU \right) } + a \int_{\Ov{\Omega}} 
	\D \ {\rm trace}[\mathcal{R}]    \right]  +\ \Lambda_\vr \left[ \intO{ E \left(\vr, \vu \ \Big|\  r, \vU \right) } + a \int_{\Ov{\Omega}} 
	\D \ {\rm trace}[\mathcal{R} ]   \right] \br
	&  + \nu \| \vu - \vU \|^2_{W^{1,2}_0(\Omega; R^3)} +  \Lambda_\vu \|\vu - \vU \|^2_{L^2(\Omega; R^3)} +  \frac{\Lambda_\vu}{4} 
	\intO{ \vr |\vu - \vU|^2 }
	\br 
	&\quad \leq \Big[ \Gamma (\| {\rm data} \|) + \chi \Big] \left[ \intO{ E \left(\vr, \vu \ \Big|\  r, \vU \right) } + a \int_{\Ov{\Omega}} 
	\D \ {\rm trace}[\mathcal{R}]    \right] + 1\br
	&\quad \quad +  \Lambda_\vr \intO{ (\vr - r)  \vU  
		\cdot(\vu - \vU) } \br&\quad \quad + \Lambda_\vr \intO{ (I_\delta[r] - r)  \vU  
		\cdot(\vU - \vu) } \br &\quad \quad+   \Lambda_\vr \intO{ \Big(P(I_\delta[r]) -  P'(r) (I_\delta[r] - r) - P(r) \Big)  } \br  &\quad \quad+  \Lambda_\vu \intO{ (1 + \vr) (I_\delta[\vU] - \vU) \cdot (\vu - \vU) }, 
	\label{b9}
\end{align} 
where 
\begin{equation} \label{b9a}
	\chi \geq 0,\ \| \chi \|_{L^1(0,T)} \leq \Gamma (\| {\rm data} \| ).
\end{equation}

Finally, similarly to the above, we have
\begin{align} 
 \frac{1}{2}  & \Lambda_\vr \intO{ (r - \vr)   \vU 
	\cdot(\vU - \vu) } = 	 \frac{1}{2}  \Lambda_\vr \int_{\vr \leq 2r} (r - \vr)   \vU 
	\cdot(\vU - \vu) \ \dx  \br &+ 
 \frac{1}{2}  \Lambda_\vr \int_{\vr > 2r} (r - \vr)   \vU 
\cdot(\vU - \vu) \ \dx,
\nonumber
\end{align}	
where 
\begin{align}
\frac{1}{2}  \Lambda_\vr \int_{\vr \leq 2r} (r - \vr)   \vU 
\cdot(\vU - \vu) \ \dx \leq \omega \Ov{U}^2 \Lambda_\vr \int_{\vr \leq 2r} 
(r - \vr)^2 \ \dx + c(\omega) \Lambda_\vr \| \vu - \vU \|^2_{L^2(\Omega; R^d)} 		
\label{b10}
\end{align}	
for any $\omega > 0$. As for the second integral, we have 
\begin{align}
\frac{1}{2}  &\Lambda_\vr \int_{\vr > 2r} (r - \vr)   \vU 
\cdot(\vU - \vu) \ \dx \leq 
\frac{1}{2} \Lambda_\vr \Ov{U} \int_{\vr > 2r} \sqrt{\vr} \sqrt{\vr} 
|\vu - \vU| \ \dx \br 
&\leq \omega \Lambda_\vr \Ov{U}^2 \int_{\vr > 2r } \vr \ \dx + 
c(\omega) \Lambda_\vr  \intO{ \vr |\vu - \vU|^2 } 
\label{b11}
\end{align}
for any $\omega > 0$. Thus choosing $\omega > 0$ small enough and 
\begin{equation} \label{b12}
	\Lambda_\vu \geq \Gamma(\| {\rm data} \| ) \Lambda_\vr
\end{equation} 	
we may write \eqref{b9} in the final form 
\begin{align} 
 &\frac{\D }{\dt} \left[ \intO{ E \left(\vr, \vu \ \Big|\  r, \vU \right) } + a \int_{\Ov{\Omega}} 
 \D \ {\rm trace}[\mathcal{R}]    \right]  +\ \Lambda_\vr \left[ \intO{ E \left(\vr, \vu \ \Big|\  r, \vU \right) } + a \int_{\Ov{\Omega}} 
 \D \ {\rm trace}[\mathcal{R} ]   \right] \br
	&  + \nu \| \vu - \vU \|^2_{W^{1,2}_0(\Omega; R^3)} +  \frac{\Lambda_\vu}{8} \left(  \|\vu - \vU \|^2_{L^2(\Omega; R^3)} +   
	\intO{ \vr |\vu - \vU|^2 } \right)
	\br 
	&\quad \leq \Big[ \Gamma (\| {\rm data} \|) + \chi \Big]  \left[ \intO{ E \left(\vr, \vu \ \Big|\  r, \vU \right) } + a \int_{\Ov{\Omega}} 
	\D \ {\rm trace}[\mathcal{R} ]   \right] + 1\br
	&\quad \quad + \Lambda_\vr \intO{ (I_\delta[r] - r)  \vU  
		\cdot(\vU - \vu) } \br &\quad \quad+   \Lambda_\vr \intO{ \Big(P(I_\delta[r]) -  P'(r) (I_\delta[r] - r) - P(r) \Big)  } \br  &\quad \quad+  \Lambda_\vu \intO{ (1 + \vr) (I_\delta[\vU] - \vU) \cdot (\vu - \vU) }.
	\label{b13}
\end{align}

\subsection{Bounds on the interpolation error}

Keeping in mind \eqref{b12} we may absorb partially the three last integrals in \eqref{b13} by 
the dissipation term
\[
\nu \| \vu - \vU \|^2_{W^{1,2}_0(\Omega; R^3)} +  \frac{\Lambda_\vu}{8} \left(  \|\vu - \vU \|^2_{L^2(\Omega; R^3)} +   
\intO{ \vr |\vu - \vU|^2 } \right)
\]
Consequently, inequality \eqref{b13} gives rise to   
\begin{align} 
	&\frac{\D }{\dt} \left[ \intO{ E \left(\vr, \vu \ \Big|\  r, \vU \right) } + a \int_{\Ov{\Omega}} 
	\D \ {\rm trace}[\mathcal{R}]    \right] \br& +  \frac{\Lambda_\vr}{2} \left[ \intO{ E \left(\vr, \vu \ \Big|\  r, \vU \right) } + a \int_{\Ov{\Omega}} 
	\D \ {\rm trace}[\mathcal{R}]    \right] \br 
	&\quad \leq \Big[ \Gamma (\| {\rm data} \|) + \chi \Big] \left[ \intO{ E \left(\vr, \vu \ \Big|\  r, \vU \right) } + a \int_{\Ov{\Omega}} 
	\D \ {\rm trace}[\mathcal{R}]    \right] + 1\br
	&\quad \quad + \Lambda_\vr \| I_\delta[r] - r \|^2_{L^\infty((0,T) \times \Omega)}  + 
	\Lambda_\vu \| I_\delta[\vU] - \vU \|^2_{L^\infty((0,T) \times \Omega; R^3)}.
	\label{b13a}
\end{align}

Now, as $q > 8$, the Sobolev space $W^{1,q}((0,T) \times \Omega)$ is embedded in the space of H\" older continuous function 
$C^{\beta} ([0,T] \times \Ov{\Omega})$, where $\beta = 1 - \frac{4}{q} > \frac{1}{2}$. Consequently, a short inspection of the 
definition on the interpolation operators in Section \ref{d} reveals 
\begin{align} 
\| I_\delta[r] - r \|^2_{L^\infty((0,T) \times \Omega)} &\leq \delta \| r \|^2_{C^\beta([0,T] \times \Ov{\Omega})} \leq 
\delta \Gamma (\| {\rm data} \|) \br 
\| I_\delta[\vU] - \vU \|^2_{L^\infty((0,T) \times \Omega)} &\leq \delta \| \vU \|^2_{C^\beta([0,T] \times \Ov{\Omega} ; R^3)} \leq 
\delta \Gamma (\| {\rm data} \|) .
\label{b13b} 
\end{align}

Choosing $\delta > 0$ small enough we get
\begin{equation} \label{b14}
\delta \Lambda_\vu  \Gamma (\| {\rm data} \|) \leq 1 \ \Rightarrow \ 
\delta \Lambda_\vr  \Gamma (\| {\rm data} \|) \leq 1. 
\end{equation}
Consequently, inequality \eqref{b13} yields
\begin{align} 
\frac{\D }{\dt} &\left[ \intO{ E \left(\vr, \vu \ \Big|\  r, \vU \right) } + a \int_{\Ov{\Omega}} 
\D \ {\rm trace}[\mathcal{R}]    \right]  + \frac{\Lambda_\vr}{2} \left[ \intO{ E \left(\vr, \vu \ \Big|\  r, \vU \right) } + a \int_{\Ov{\Omega}} 
\D \ {\rm trace}[\mathcal{R} ]   \right] \br
&\leq \chi \left[ \intO{ E \left(\vr, \vu \ \Big|\  r, \vU \right) } + a \int_{\Ov{\Omega}} 
\D \ {\rm trace}[\mathcal{R}]    \right] + 3, 
\label{b15} 
\end{align}
where we have identified $\chi \approx \Gamma (\| {\rm data} \|) + \chi$
still satisfying \eqref{b9a}. 
Thus a direct manipulation yields
\begin{align}
\exp &\left( \int_0^\tau \left( \frac{\Lambda_\vr}{2} - \chi(s) \right) \D s \right) \left[ \intO{ E \left(\vr, \vu \ \Big|\  r, \vU \right) (\tau, \cdot) } + a \int_{\Ov{\Omega}} 
\D \ {\rm trace}[\mathcal{R}](\tau, \cdot)    \right] \br
&\leq 3 \int_0^\tau \exp \left( \int_0^t \left( \frac{\Lambda_\vr}{2} - \chi(s) \right) \D s \right) \dt  \br  &+
\intO{ E \left(\vr, \vu \ \Big|\  r, \vU \right)(0, \cdot)} \leq \frac{6}{\Lambda_\vr} \exp \left( \frac{\Lambda_\vr}{2} \tau \right) 
+ \intO{ E \left(\vr, \vu \ \Big|\  r, \vU \right)(0, \cdot)},
\nonumber 
\end{align}
and, consequently, 
\begin{align} \label{b16}
&\left[ \intO{ E \left(\vr, \vu \ \Big|\  r, \vU \right) (\tau, \cdot) } + a \int_{\Ov{\Omega}} 
\D \ {\rm trace}[\mathcal{R}](\tau, \cdot)    \right]	\br
&\quad \leq \frac{6}{\Lambda_\vr} \exp \left( \int_0^\tau \chi(s) \D s \right)
+ \exp \left( - \frac{\Lambda_\vr}{2} \tau \right) \exp \left( \int_0^\tau \chi(s) \D s \right) \intO{ E \left(\vr, \vu \ \Big|\  r, \vU \right)(0, \cdot)} \br 
&\leq \left[ \frac{2}{\Lambda_\vr} + \exp \left( - \frac{\Lambda_\vr}{2} \tau \right) \right] \Gamma (\| {\rm data} \|)
\end{align}	
for a.a. $0 \leq \tau \leq T$.

\section{Control in the prediction (forecast) period}
\label{C}

Finally, we estimate the distance of $(r, \vU)$ and $(\vr, \vu)$ in the prediction (forecast) period $(T,T^+)$.
We focus on the more difficult and physically relevant case $\gamma \leq 2$.
As the nudging is no longer active, we have the standard form of the relative energy inequality:
\begin{align} 
	&\frac{\D }{\dt} \left[ \intO{ E \left(\vr, \vu \ \Big|\  r, \vU \right)  } + a \int_{\Ov{\Omega}} 
	\D \ {\rm trace}[\mathcal{R}]   \right] + \nu \| \vu - \vU \|^2_{W^{1,2}_0(\Omega; R^d)}  
\br 
	&\quad \leq - \intTd{ \vr (\vU - \vu) \cdot \Grad \vU \cdot (\vU - \vu) } \br
	&\quad \quad -  \intO{ \Div \vU \Big( p(\vr) - p'(r) (\vr - r) -  p(r) \Big) } \br
	&\quad \quad + \int_{\Omega} \Grad \vU : \D \mathcal{R} 
	 \br
	&\quad \quad + \intO{ (\vr - r) \Big( \frac{1}{r} \Div \mathbb{S}(\Ds \vU) + \vc{g} \Big) \cdot (\vU - \vu) }                  \br
	&\quad \leq \| \vU \|_{L^\infty(T,T^+; W^{1, \infty}(\Omega; R^d))} \left[ \intO{ E \left(\vr, \vu \ \Big|\  r, \vU \right)  } + a \int_{\Ov{\Omega}} 
	\D \ {\rm trace}[\mathcal{R}]   \right]
	 \br 
&\quad \quad + 	c(\nu, \Ov{r}) \left\| \frac{1}{r} \Div \mathbb{S}(\Ds \vU) + \vc{g} \right\|_{L^3(\Omega)}^2 \intO{ E \left(\vr, \vu \Big| r, \vU \right) } + \frac{\nu}{2} \| \vu - \vU \|^2_{W^{1,2}_0(\Omega; R^d)} \br 
&\quad \quad + \int_{\vr > 2r} (\vr - r) \Big( \frac{1}{r} \Div \mathbb{S}(\Ds \vU) + \vc{g} \Big) \cdot (\vU - \vu) \dx.
 	\label{C1}
\end{align} 
The last integral is handled by means of H\" older inequality for 
$p = \gamma, q = \frac {6\gamma} {5 \gamma - 6}, r =6$:
\begin{align}
\int_{\vr > 2r} &(\vr - r) \Big( \frac{1}{r} \Div \mathbb{S}(\Ds \vU) + \vc{g} \Big) \cdot (\vU - \vu) \dx \br
&\leq \left( \int_{\vr > 2r} |\vr - r|^\gamma \ \dx \right)^{\frac{1}{\gamma}} 
 \left\| \frac{1}{r} \Div \mathbb{S}(\Ds \vU) + \vc{g} \right\|_{L^q(\Omega; R^{3 \times 3})} \| \vU - \vu \|_{L^6(\Omega; R^3)}. 
\label{C1a}
\end{align}

\begin{Remark} \label{IR2}
	
This is the only part of the proof, where the hypothesis $\gamma > \frac{6}{5}$ is used. This restriction could be relaxed if 
the observed solution was more regular, specifically 
\[
\Div \mathbb{S}(\Ds \vU) \in L^2(T^-, T^+; L^\infty(\Omega; R^3)). 
\]	
This is indeed true if 
\[
r(T^-, \cdot) \in W^{3,2}(\Omega), \ \vU (T^-, \cdot) \in W^{3,2}(\Omega; R^3) \ \mbox{+ compatibility conditions},
\]
see \cite{BaFeMi}.

\end{Remark}	

Combining \eqref{C1} with \eqref{C1a} we obtain 
\begin{align} 
	\frac{\D }{\dt} &\left[ \intO{ E \left(\vr, \vu \ \Big|\  r, \vU \right)  } + a \int_{\Ov{\Omega}} 
	\D \ {\rm trace}[\mathcal{R}]   \right] \br & \leq \chi(t) \left( \left[ \intO{ E \left(\vr, \vu \ \Big|\  r, \vU \right)  } + a \int_{\Ov{\Omega}} 
	\D \ {\rm trace}[\mathcal{R}]   \right] 
		+  \left( \intO{ E \left(\vr, \vu \Big| r, \vU \right) } \right)^{\frac{2}{\gamma}} \right)
	\label{C3}
\end{align} 
for any $T \leq t \leq T^+$, where $\chi$ satisfies \eqref{b9a}. Thus we get the estimate 
\begin{align} 
 &\left[ \intO{ E \left(\vr, \vu \ \Big|\  r, \vU \right) (\tau, \cdot)  } + a \int_{\Ov{\Omega}} 
 \D \ {\rm trace}[\mathcal{R}](\tau, \cdot)   \right] \br
&\quad \leq \exp\left(2 \int_{T}^\tau \chi(s) \D s    \right)\left[ \intO{ E \left(\vr, \vu \ \Big|\  r, \vU \right) (T, \cdot)  } + a \int_{\Ov{\Omega}} 
\D \ {\rm trace}[\mathcal{R}](T, \cdot)   \right] 
\label{C4}
\end{align}
for $T < \tau \leq T^+$ as long as 
\begin{equation} \label{C5}
	\exp\left(2 \int_{T}^{T^+} \chi(s) \D s    \right)\left[ \intO{ E \left(\vr, \vu \ \Big|\  r, \vU \right) (T, \cdot)  } + a \int_{\Ov{\Omega}} 
	\D \ {\rm trace}[\mathcal{R}](T, \cdot)   \right]  \leq 1.
\end{equation}	

However, the value of the integral 
\[
\left[ \intO{ E \left(\vr, \vu \ \Big|\  r, \vU \right) (T, \cdot)  } + a \int_{\Ov{\Omega}} 
\D \ {\rm trace}[\mathcal{R}](T, \cdot)   \right] 
\]
is controlled by \eqref{b16}, specifically, 
\begin{align}
&	\left[ \intO{ E \left(\vr, \vu \ \Big|\  r, \vU \right) (\tau, \cdot)  } + a \int_{\Ov{\Omega}} 
	\D \ {\rm trace}[\mathcal{R}](\tau, \cdot)   \right] \br
& \quad \leq \exp\left(2 \int_{T}^\tau \chi(s) \D s    \right)  \left[ \frac{2}{\Lambda_\vr} + \exp \left( - \frac{\Lambda_\vr}{2} T \right) \right] \Gamma (\| {\rm data} \|)  
\approx \left[ \frac{2}{\Lambda_\vr} + \exp \left( - \frac{\Lambda_\vr}{2} T \right) \right] \Gamma (\| {\rm data} \|) 
\nonumber
\end{align}	
for a.a. $0\leq  \tau \leq T$. Replacing $\Lambda_\vr$ by $2 \Lambda_\vr$ as the case may be, we have completed the proof of Theorem \ref{mT1}.


\def\cprime{$'$} \def\ocirc#1{\ifmmode\setbox0=\hbox{$#1$}\dimen0=\ht0
	\advance\dimen0 by1pt\rlap{\hbox to\wd0{\hss\raise\dimen0
			\hbox{\hskip.2em$\scriptscriptstyle\circ$}\hss}}#1\else {\accent"17 #1}\fi}

\end{document}